\documentclass{article}
\usepackage{amsmath}
\usepackage{cite}
\newtheorem{Thm}{Theorem}

\newcommand{\Dt}{|S(t+h) - S(t)|}

\begin{document}
\title{Gram's Law Fails a Positive Proportion of the Time}
\author{Timothy Trudgian\thanks{The author is gracious of the support of a General Sir John Monash Award.}\\
Mathematical Institute\\
University of Oxford, United Kingdom}
\maketitle

\begin{abstract}

It is straightforward to show that all the non-trivial zeroes of the Riemann zeta-function $\zeta(s)$ are confined to the critical strip: $0\leq \Re(s) \leq 1.$ It is another matter to seek out their precise location. Extending the work of Riemann, Gram introduced a procedure for detecting these zeroes: Gram's Law. It is known that Gram's Law fails infinitely often, and that a weaker formulation of Gram's Law is true infinitely often. This paper extends these results by showing that there is a positive proportion of both failures and (weak) successes.
\end{abstract}
\section{Gram's Law}

Connected to $\zeta(s)$ are two functions known as the Riemann-Siegel functions,
\begin{equation}\label{RS}
Z(t) = e^{i\theta(t)} \zeta\left(\frac{1}{2} +it\right),
\end{equation}
where
\begin{equation}\label{the2}
\theta(t) \sim \frac{t}{2}\log\frac{t}{2\pi} + O(t),
\end{equation}
is a steadily increasing function of $t$ and $Z(t)$ is real-valued whenever $t$ itself is real-valued. Hence the zeroes of $\zeta(\frac{1}{2} +it)$ coincide precisely with those of $Z(t)$. Therefore the search for zeroes of the zeta-function is equivalent to a search for changes in sign of $Z(t)$. From (\ref{RS}) it follows that
\begin{equation}
\zeta\left(\frac{1}{2}+it\right) = Z(t)\cos\{\theta(t)\} -iZ(t)\sin\{\theta(t)\},
\end{equation}
and so at the \textbf{Gram points} $\theta(g_{n}) =n\pi$
\begin{equation}
\zeta\left(\frac{1}{2} +ig_{n}\right) =(-1)^{n} Z(g_{n}).
\end{equation}
Gram \cite{Gram} found that $\Re\{\zeta(\frac{1}{2} +it)\} >0$ frequently\footnote{The Riemann-Siegel formula gives some basis for this empirical observation: $(-1)^{n}Z(g_{n}) = 2\sum_{\nu=1}^{X} \nu^{-1/2}\cos(\nu\log g_{n}) + O(g_{n}^{-1/4})$, where $X=(g_{n}/2\pi)^{1/2}$. The first term in the summation is $+1$ and thereafter the terms are oscillatory and decreasing in magnitude.} so that, in particular $Z(g_{n})$ and $Z(g_{n+1})$ were often of opposite sign, whence a zero of $\zeta(\frac{1}{2} +it)$ must occur in this interval. Gram noted that this pattern can reasonably be expected to continue: one zero of $\zeta(\frac{1}{2} +it)$ between successive Gram points.
Gram located the first $15$ non-trivial zeroes of $\zeta(s)$ using this method and found each one to lie on the line. Hutchinson extended these results in \cite{Hutchinson} to show that the first 138 zeroes are on the line and later Titchmarsh \cite{Titchmarsh4} improved this to the first 1041 zeroes. In all cases it was shown that there were no other complex zeroes up to these heights. For a comprehensive historical account of the applications of Gram's Law to finding zeroes of the zeta-function, see \cite[pp.171-182]{Edwards}.

\textbf{Gram's Law}, as defined by Hutchinson (\textit{ibid}) is that statement that there is \textit{exactly} one zero of $\zeta(\frac{1}{2} +it)$ for $t\in(g_{n}, g_{n+1}]$ and that this is the only\footnote{Thus Gram's Law if it were true implies the Riemann hypothesis.} zero of with $\Re(t) \in (g_{n}, g_{n+1}]$. A weakened version of Gram's Law (hereafter referred to as the \textbf{Weak Gram Law}) is the statement that there is \textit{at least} one zero of $\zeta(\frac{1}{2} +it)$ for $t\in(g_{n}, g_{n+1}]$. The successes of Gram's Law were first investigated by analysing the discrete properties of the function $Z(t)$, but later the argument function $S(t)$ was used to obtain improvements.

\subsection{The function $S(t)$ and Gram's Law}

As is standard let $S(T) = \pi^{-1}\arg\{\zeta(\frac{1}{2} +iT)\}$ (for more properties on the function $S(T)$ see, for example \cite[pp. 212-223]{Titchmarsh}) and so by the Riemann-von Mangoldt formula
\begin{equation}
S(T) = N(T) -\pi^{-1}\theta(T) -1,
\end{equation}
where $N(T)$ is the number of zeroes of $\zeta(\sigma +it)$ for $0\leq t\leq T$. Since $\theta(g_{n}) = n\pi$ and $\theta(t)$ is a steadily increasing function of $t$, it is seen at once that $S(T)$ is integral precisely at the Gram points. Furthermore if $S(g_{n}) = \lambda$ and there is exactly one zero in the Gram interval $(g_{n}, g_{n+1}]$ then $S(g_{n+1}) = \lambda$ and $|S(t) - \lambda|\leq 1$ throughout the interval. So intervals in which Gram's Law is valid induce some constancy in the function $S(t)$ and it is this constancy which forms the basis of the following analysis.

\section{General Failures}\label{s1}

Titchmarsh showed in \cite{Titchmarsh1} that Gram's Law fails infinitely often. This section will show that Gram's Law and the Weak Gram Law fail a positive proportion of the time (given in \textbf{Theorem \ref{t1}} on p.\pageref{t1} and \textbf{Theorem \ref{t2}} on p.\pageref{t2} respectively).

The following result due to Fujii \cite{Fujii} concerns the `shifted moments' of $S(t)$, viz. 
\begin{align}\label{e1}
I(T) &= \int_{T}^{2T} \Dt^{2} \, dt \notag \\ 
&= \pi^{-2} T\log\left(3+h\log T\right) + O\left[T\{\log\left(3+ h\log T\right)\}^{\frac{1}{2}}\right],
\end{align}
which is valid for $0 \leq h\leq \frac{1}{2}T$. This becomes an asymptotic relationship, i.e.
\begin{equation}
I(T) \sim\pi^{-2} T\log\left(3+h\log T\right),
\end{equation}
if $h\log T$ is sufficiently large. Henceforth $h=C_{0}(\log T)^{-1}$, where $C_{0}$ is a constant that is chosen to be sufficiently large to ensure the dominance of the main term in (\ref{e1}) over the error term.

If $t$ and $t+h$ are in a connected union of Gram intervals in which Gram's Law is valid, then $|S(t+h)-S(t)| \leq 2$. Thence $I(T) \leq 4T$ which is `too small', in that this is not asymptotic to $\pi^{-2} T\log\left(3+h\log T\right)$. This lends credence to what has already been shown using the work of Ghosh \cite{Ghosh}: that for sufficiently large $T$ there must be at least one failure\footnote{This has been shown by the author and is being prepared for publication.} of Gram's Law between heights $T$ and $2T$. In some loose sense, if $S$ is the set on which Gram's Law is valid and $\overline{S}$ is the complement of $S$ in $[T, 2T]$, then (\ref{e1}) can be rewritten
\begin{align}\label{int}
I(T) &\sim\pi^{-2} T\log\left(3+h\log T\right)\notag \\
&= \int_{S} \Dt^{2}\, dt + \int_{\overline{S}} \Dt^{2} \, dt \notag \\
&\leq 4|S| + \int_{\overline{S}} \Dt^{2} \, dt,
\end{align}
whence an estimate on $|\overline{S}|$ can be made.

To this end, let the sequences $\{i_{n}\}$ and $\{j_{n}\}$ index the Gram points such that Gram's Law holds on the collection of intervals $\left(g_{i_{n}}, g_{j_{n}}\right]$ and Gram's Law fails on the collection of intervals $\left(g_{j_{n}}, g_{i_{n+1}}\right]$. Also let $k_{n} = i_{n+1} - j_{n}$, that is, the number of consecutive Gram points between which Gram's Law fails. So then $ \sum_{n} k_{n} = N_{F}(2T)$: the number of failures between heights $T$ and $2T$. It is now appropriate to introduce the following elementary result concerning $N_{g}$: the number of Gram points between heights $T$ and $2T$,
\newtheorem{lem}{Lemma}\label{firstl}
\begin{lem}
$N_{g}\sim (2\pi)^{-1}T\log T$. Furthermore if $g_{n}$ and $g_{m}$ are Gram points in the interval $[T, 2T]$ then $g_{n}-g_{m} = O(\frac{n-m}{\log T})$.
\end{lem}
The first statement follows from both (\ref{the2}) and the definition of the Gram points. Since it can be shown (see \cite[p 263]{Titchmarsh}) that $\theta'(t) \sim\frac{1}{2}\log t$, then the mean value theorem gives
\begin{equation}
 \frac{\theta(g_{n}) - \theta(g_{m})}{g_{n}-g_{m}} = \frac{(n-m)\pi}{g_{n}-g_{m}} = \frac{1}{2}\log\xi,
\end{equation}
for some $\xi\in(g_{n}, g_{m})$, whence the result follows.

It is clear that the relative locations of $t$ and $t+h$ will determine the bound on $\Dt$: namely if $g_{i_{n}} \leq t\leq t+h\leq g_{j_{n}}$ then $\Dt \leq 2$. This leads to the definition
\begin{equation}
S:= \{t\in[T, 2T] \quad : \exists n \quad : g_{i_{n}} \leq t\leq t+h\leq g_{j_{n}}\},
\end{equation}
whence $\int_{S} \Dt^{2}\, dt \leq 4T$, as claimed. 

Now let $\overline{S}$ be the complement of $S$ in $[T, 2T]$. Then, if $t$ belongs to $\overline{S}$ either $t\in \left(g_{i_{n}}, g_{j_{n}}\right]$ and $t+h \geq g_{j_{n}}$; or $t\in \left(g_{j_{n}}, g_{i_{n+1}}\right]$. The former condition is equivalent to $g_{j_{n}} \geq t\geq g_{j_{n}} - h$ and so in any case $g_{j_{n}} - h\leq t\leq g_{i_{n+1}}.$ These intervals may overlap in $[T, 2T]$ and indeed 
\begin{equation}
\overline{S} \subset \bigcup_{n} \left(g_{j_{n}} - h, g_{i_{n+1}}\right].
\end{equation}
Whether or not these intervals are disjoint is of no consequence for \textbf{Lemma 1} gives
\begin{equation}\label{uS}
|\overline{S}| \ll \sum_{n} h+ \frac{k_{n}}{\log T} \ll \left(h+ \frac{1}{\log T}\right)N_{F}(2T).
\end{equation}

Ultimately an estimate on this number $N_{F}(2T)$ is sought and hence the imposition of a lower bound of (\ref{uS}) would be useful. Returning to (\ref{int}) it is seen that
\begin{align}\label{LHS}
\pi^{-2} T\log\left(3+h\log T\right) &\leq 4|S| + \int_{\overline{S}} \Dt^{2} \, dt \notag \\
&\leq 4T + \int_{\overline{S}} \Dt^{2} \, dt.
\end{align}
Currently $h=C_{0}(\log T)^{-1}$ and $C_{0}$ is chosen to be sufficiently large such that the main term in (\ref{e1}) dominates the error term. If, in addition to this, $C_{0}$ is taken large enough to make the quantity $\pi^{-2}T\log(3+h\log T)$ large than $5T$, then (\ref{LHS}) gives
\begin{equation}\label{8}
T\ll \int_{\overline{S}} \Dt^{2} \, dt.
\end{equation}

Results on higher moments of the function $S(t)$ have been developed by Fujii following the work of Selberg, as detailed in \cite[pp. 245-246]{Titchmarsh}. These may be employed after an application of Cauchy's inequality to given
\begin{align}\label{cd1}
\int_{\overline{S}} \Dt^{2} \, dt &\leq \left(\int_{\overline{S}} \Dt^{4} \, dt\right)^{\frac{1}{2}}\times\left(\int_{\overline{S}}\, dt\right)^{\frac{1}{2}}\notag \\
&= |\overline{S}|^{\frac{1}{2}}\times\left(\int_{\overline{S}} \Dt^{4} \, dt\right)^{\frac{1}{2}},
\end{align}
whence via (\ref{8}) it follows that
\begin{equation}\label{rex}
 T\ll |\overline{S}|^{\frac{1}{2}}\times\left(\int_{\overline{S}} \Dt^{4} \, dt\right)^{\frac{1}{2}}.
\end{equation}
The particular result needed here is
\begin{equation}\label{cd2}
\int_{T}^{2T} |S(t+h) - S(t)|^{4}\, dt \ll T \log^{2}\left(3+h\log T\right).
\end{equation}
Since $\int_{\overline{S}} |\cdot | \leq \int_{T}^{2T} |\cdot |$, the inequality in (\ref{rex}) can be combined with (\ref{cd2}) to give
\begin{equation}
T \ll |\overline{S}|^{\frac{1}{2}} T^{\frac{1}{2}} \log\left(3+h\log T\right),
\end{equation}
and so a lower bound on $|\overline{S}|$ is attained, viz.
\begin{equation}\label{lS}
|\overline{S}| \gg \frac{T}{\log^{2}\left(3+h\log T\right)}.
\end{equation}
Now the upper bound for $|\overline{S}|$ in (\ref{uS}) can be combined with the lower bound in (\ref{lS}) to give,
\begin{equation}\label{13}
\frac{T}{\log^{2}\left(3+h\log T\right)} \ll \left(\frac{1}{\log T} + h\right)N_{F}(2T),
\end{equation}
or
\begin{equation}\label{14}
N_{F}(2T) \gg \frac{T\log T}{(1+C_{0})\log^{2}(3+C_{0})} \gg AT\log T,
\end{equation}
for some positive constant $A$. By \textbf{Lemma 1} the total number of Gram points between heights $T$ and $2T$ is $O(T\log T)$ which proves the following
\begin{Thm}\label{t1}
For sufficiently large $T$ there is a positive proportion of failures of Gram's Law between $T$ and $2T$. 
\end{Thm}

\subsection{Further failures}
A small point to note is that it is now possible to deduce that there is a positive proportion of Gram intervals which \textbf{do not} contain a zero of $\zeta(\frac{1}{2} +it)$. For $i=0, 1, 2, \ldots$ let $F_{i}$ denote a Gram interval in which $i$ zeroes are located. Furthermore, let $N_{F_{i}}$ denote the number of Gram intervals between heights $T$ and $2T$ which contain exactly $i$ zeroes: so that $N_{F_{0}}$ is the number of $F_{0}$ intervals, $N_{F_{1}}$ the number of intervals in which Gram's Law is valid, and so on. Then
\begin{equation}\label{i}
N_{F_{0}} + N_{F_{1}} + \ldots + N_{F_{k}} + \ldots = N_{g} = \frac{T}{2\pi}\log T+ O(T) ,
\end{equation}
where $N_{g}$ is the total number of Gram intervals between heights $T$ and $2T$. That there is a positive proportion of failures is represented by the following equation
\begin{equation}\label{ii}
N_{F_{0}} + N_{F_{2}} + \ldots + N_{F_{k}} + \ldots \geq AN_{g},
\end{equation}
where, as before, the number $N_{F_{1}}$ is absent since this does not represent any failures. Lastly since all the zeroes on the critical line between heights $T$ and $2T$, denoted by $N_{0} (T)$, fall within Gram intervals a third relation may be written, viz.
\begin{equation}\label{iii}
N_{F_{1}} + 2N_{F_{2}} + \ldots + kN_{F_{k}} +\ldots = N_{0} (T) \leq N(T) = \frac{T}{2\pi}\log T + O(T),
\end{equation}
where $N(T)$ is the number of complex zeroes of $\zeta(s +it)$ with $0\leq t\leq T$.
The subtraction of equation (\ref{i}) from (\ref{iii}) gives
\begin{align}
O(T) &\geq -N_{F_{0}} + N_{F_{2}} + 2N_{F_{3}} \ldots + (k-1)N_{F_{k}} + \ldots\\
&\geq -N_{F_{0}} + N_{F_{2}} + N_{F_{3}} \ldots +N_{F_{k}} + \ldots
\end{align}
whence, upon an addition of $2N_{F_{0}}$ and an invocation of (\ref{ii}) it is seen that
\begin{equation}
2N_{F_{0}} + O(T) \geq N_{F_{0}} + N_{F_{2}} + N_{F_{3}} + \ldots N_{F_{k}} + \ldots \geq AN_{g},
\end{equation}
so that
\begin{equation}
\frac{N_{F_{0}}}{N_{g}} \geq \frac{A}{2} + O\left(\frac{1}{\log T}\right).
\end{equation}
Thus the following has been now been proved
\begin{Thm}\label{t2}
For sufficiently large $T$ there is a positive proportion of failures of the Weak Gram Law between $T$ and $2T$. 
\end{Thm}

Since the number of $F_{0}$ intervals is certainly less than the total number of violations of Gram's Law, the order of $N_{F_{0}}$ is exactly determined, viz. $AT\log T \leq N_{F_{0}} \leq AT\log T$. There is little else\footnote{One possibility is to calculate these constants, but this is not achievable via the methods in this paper.} to be said about the nature of $F_{0}$ intervals, so it is natural to now turn to the remaining cases: those Gram intervals which contain at least one zero of $\zeta(\frac{1}{2} +it)$.

\section{$F_{k}$ intervals}

Titchmarsh showed in \cite{Titchmarsh2} that the Weak Gram Law is true\footnote{One remark to be made is that it is not yet known whether Gram's Law is true infinitely often.} infinitely often. What is actually shown in his proof is that there is an infinite number of Gram intervals which contain an \textbf{odd number} of zeroes. His proof concludes that the proportion of Gram intervals between $T$ and $2T$ which contain an odd number of zeroes of $\zeta(\frac{1}{2} +it)$ is greater than $A(T^{1/3}\log^{2}T)^{-1}$, with $A$ a positive constant. This section will show (in \textbf{Theorem \ref{the3}} on p.\pageref{the3}) that the Weak Gram Law is true a positive proportion of the time.

\subsection{Outline}

It is difficult to investigate the quantities $N_{F_{k}}$ for `small' $k$, since the induced behaviour in $S(t)$ is virtually undetectable. Indeed the methods used in \S 1 viz. shifted moments of $S(t)$ are unable to distinguish a collection of $F_{1}$ intervals from a sequence of alternating $F_{0}$ and $F_{2}$ intervals. Investigations into the frequency of successes of Gram's Law (or the quantity $N_{F_{1}}$) must be made through some other route. What can be said is a measure of the success of Gram's Law in its weak sense: that is the number of intervals which contain \textit{at least} one zero of $\zeta(\frac{1}{2} +it)$.

By Selberg's result, a positive proportion of zeroes lie on the critical line; clearly each zero is contained within a Gram interval. There is a possibility that when $k$ is arbitrarily large there are many $F_{k}$ intervals which contain the bulk of these zeroes. To rule out this possibility it is necessary to place a bound on the growth of $N_{F_{k}}$ as $k\rightarrow \infty$. Once this has been established it will be shown that there is a $K$ such that the Gram intervals containing fewer than $K$ zeroes together contain the positive proportion of zeroes.

\subsection{Improvements in the function $S(t)$}

Much work has been done concerning the number of zeroes of $\zeta(\frac{1}{2} +it)$ of multiplicity greater than one. Extending this work to short intervals, particularly Gram intervals is natural since a zero of order $m$ will induce similar behaviour in $S(t)$ as will $m$ simple zeroes. The following result is due to Korolev \cite{Korolev2005}
\begin{equation}\label{evenm}
I_{m} = \int_{T}^{T+H}|S(t+h) - S(t)|^{2m} \, dt \leq \left(Cm^{2}\right)^{m}H,
\end{equation}
where
\begin{equation}
H = T^{\frac{27}{82}} +\epsilon; \qquad 0<\epsilon < 0.001; \qquad h = \frac{2\pi}{3\log{\frac{T}{2\pi}}},
\end{equation}
and $C$ is given as an explicit positive constant. This formula has $h=c\,(\log T)^{-1}$, where $c$ is small relative to the length of Gram intervals. In order to easily detect the contribution of an $F_{k}$ interval to the integrand in (\ref{evenm}), the $h$ must be replaced with $Nh$ (with $N$ to be chosen later) such that $Nh$ is longer than a Gram interval. The following, which is easily deduced from  \textbf{Lemma 1}, will prove useful
\begin{lem}\label{prv}
 Denote the length of the longest Gram interval in $[T, T+H]$ by $L^{+}$ and the length of the shortest by $L^{-}$. Then
\begin{equation}
 L^{+} = \frac{2\pi}{\log\frac{T}{2\pi}},
\end{equation}
and
\begin{equation}
 L^{-} =\frac{2\pi}{\log\frac{T+H}{2\pi}} = \frac{2\pi}{\log\frac{T}{2\pi} + \log(1+\frac{H}{T})} = \frac{2\pi}{\log\frac{T}{2\pi}}\left\{1+o(T)\right\}.
\end{equation}
\end{lem}

Now suppose the interval $(g_{n}, g_{n+1}]$ is an $F_{k}$ interval and that $S(g_{n}) = \lambda$. Then $S(g_{n+1}) = \lambda +k-1$ and thenceforth $S(t)$ can decrease by at most one on the interval $(g_{n+1}, g_{n+2}]$. Furthermore for $t\in (g_{n-2}, g_{n-1}]$ it follows that $S(t) < \lambda +2$. The choice of $N$ must be made such that $t+Nh \geq g_{n+1}$, which is satisfied if
\begin{equation}
 Nh\geq(g_{n+1} - g_{n-2})\geq 3L^{+} = 9h,
\end{equation}
so that $N=9$ will suffice. When $T$ is sufficiently large, \textbf{Lemma \ref{prv}} shows $t+Nh<g_{n+2}$ and so over an interval of length $L^{-}$ the difference $|S(t+h) - S(t)|$ is now bounded below by $|k-4|$. By an application of the H\"{o}lder inequality (\ref{evenm}) then becomes
\begin{equation}
 \int_{T}^{T+H}|S(t+h) - S(t)|^{2m} \, dt \leq 9^{2m-1}\sum_{i=0}^{8}\int_{T+ih}^{T+ih+H}|S(t+h) - S(t)|^{2m}\, dt.
\end{equation}
Applying (\ref{evenm}) with $T+ih$ in place of $T$ it is seen that,
\begin{equation}
\int_{T}^{T+H}|S(t+Nh) - S(t)|^{2m} \, dt \leq (Am^{2})^{m}H\{1+o(1)\},
\end{equation}
where $A$ is a positive constant. Now suppose there are $N_{F_{k}}$ intervals between heights $T$ and $T+H$. Each one will contribute at least $|k-4|^{2m}$ in the above integrand over a length at least $L^{-}$. Thus
\begin{equation}
(Am^{2})^{m}H\{1+o(1)\} \geq \frac{\pi N_{F_{k}}(T)\, (k-4)^{2m}}{\log\frac{T}{2\pi}},
\end{equation}
or, expressed more succinctly,
\begin{equation}
\frac{N_{F_{k}}(T)}{H\log T}\ll \left(\frac{Am^{2}}{(k-4)^{2}}\right)^{m}.
\end{equation}
Now the task is to find the value of $m$ depending on $k$ that minimises the right hand inequality. Let 
\begin{equation}
F(m) = \left(\frac{Am^{2}}{(k-4)^{2}}\right)^{m},
\end{equation}
then
\begin{equation}
\frac{dF(m)}{dm} = \left(\frac{Am^{2}}{(k-4)^{2}}\right)^{m}\left\{2+ \log\left(\frac{Am^{2}}{(k-4)^{2}}\right)\right\},
\end{equation}
and clearly this stationary point $m^{*}=\frac{(k-4)}{e\sqrt{A}}$ is indeed a minimum. Since $m$ is an integer the value to be taken is whichever of $[m^{*}]$ or $[m^{*}] +1$ is the nearer to $m$. The error of such an assignment of value is $O(1)$ in the exponent and can be absorbed into the ultimate $O$-constant. Thus, given a value of $k$ the value $m=m^{*}$ gives the bound
\begin{equation}\label{res}
\frac{N_{F_{k}}(T)}{H\log T}\ll \exp(-Ak),
\end{equation}
where $A$ is a positive constant, and clearly this result remains valid if $H$ is replaced with $T$. So $F_{k}$ intervals, for large $k$ are `exponentially rare'.

Since a positive proportion of zeroes lie on the critical line, the inequality in (\ref{res}) can be used to show that a positive proportion of Gram intervals contain at least one zero. For, there is a constant $A'$ such that
\begin{equation}
0< A' <\frac{N_{F_{1}}(T) + 2N_{F_{2}}(T) + \ldots  + kN_{F_{k}}(T) + \ldots}{T\log T},
\end{equation}
and by (\ref{res}) this series on the right hand side is convergent. So, if $\delta$ is any small positive number, choose $K$ so large that the sum $(T\log T)^{-1}\sum_{k= K+1}^{\infty}kN_{F_{k}}(T)$ is not greater than $A' - \delta$. Then
\begin{equation}\label{39}
0<\delta < \frac{\sum_{k=1}^{K} kN_{F_{k}}(T)}{T\log T} < K\frac{\sum_{k=1}^{K}N_{F_{k}}(T)}{T\log T},
\end{equation}
whence the number of Gram intervals which contain at least one zero is at least $AT\log T$, with $A$ a positive constant. Thus the following has been proved
\begin{Thm}\label{the3}
For sufficiently large $T$ there is a positive proportion of successes of the Weak Gram Law between $T$ and $2T$.
\end{Thm}

\section{Concluding Remarks}
From (\ref{39}) it follows that there must be a positive proportion of at least one of the $N_{F_{k}}$'s. Intuitively one might expect $N_{F_{k}}$ to be steadily decreasing with $k$ (which would be an improvement to the estimate in (\ref{res})). If such a relation could be shown it would therefore follow that there is a positive proportion of intervals in which Gram's Law is valid.

\bibliographystyle{plain}
\bibliography{Arxiv}

\end{document}